\documentclass[11pt]{amsart}
\usepackage{xcolor}
\newlength{\basicwidth}
\newlength{\shortbasicwidth}
\newlength{\basicheight}
\numberwithin{equation}{section}

\begin{document}

\begin{center}
{{\Large{A Refinement of Vietoris' Inequality
for Cosine Polynomials}}}
\end{center}

\par\vspace*{2\baselineskip}\par
\centerline{\large  ({\em To appear in: Analysis and Applications})}

\vspace{0.6cm}
\begin{center}
HORST ALZER$^a$  \quad{and}
\quad{MAN KAM KWONG$^b$}
\footnote{The research of this author is supported by the Hong Kong Government GRF Grant PolyU 5003/12P and the Hong Kong Polytechnic University Grants G-UC22 and G-UA10}
\end{center}

\vspace{0.6cm}
\begin{center}
$^a$ Morsbacher Str. 10, 51545 Waldbr\"ol, Germany\\
\emph{e-mail:} {\tt{h.alzer@gmx.de}}
\end{center}

\vspace{0.5cm}
\begin{center}
$^b$ Department of Mathematics, The Hong Kong Polytechnic University,\\
Hunghom, Hong Kong\\
\emph{e-mail:} {\tt{mankwong@poly.edu.hk}}
\end{center}

\vspace{1.5cm}
{\bf{Abstract.}} Let
$$
T_n(x)=\sum_{k=0}^n b_k\cos(kx)
$$
with 
$$
\qquad  \qquad  b_{2k}=b_{2k+1}=\frac{1}{4^k}{2k \choose k}  \qquad 
(k\geq 0).
$$
In 1958, Vietoris proved that
$$
T_n(x)>0 \qquad{(n\geq 1; x\in (0,\pi))}.
$$
We offer the following improvement of
this result:
The inequalities
$$
T_n(x)\geq c_0 + c_1 x + c_2 x^2 >0
\qquad{(c_k\in\mathbf{R}, k=0,1,2)}
$$
hold for all  $n\geq 1$ and  $x\in (0,\pi)$ if and only if
$$
c_0=\pi^2 c_2, 
\quad{c_1=-2\pi c_2,}
\quad{0<c_2\leq \alpha},
$$
where
$$
\alpha=\min_{0\leq t<\pi} 
\frac{T_6(t)}{(t-\pi)^2}=0.12290... .
$$

\vspace{1cm}

\vspace{0.1cm}

\newpage

\section{Introduction}

In 1958, Vietoris \cite{V} published the following ``surprising and quite deep result" \cite[p. 1]{K} on inequalities  for a class of sine and cosine polynomials.

\vspace{0.3cm}
\noindent
{\bf{Proposition 1.}} \emph{If the real numbers
$a_k$ $(k=0,1,...,n)$ satisfy
$$
a_0\geq a_1 \geq \cdots \geq a_n>0
\quad{and}
\quad{2k a_{2k}\leq (2k-1) a_{2k-1} \,\,\, (k\geq 1)},
$$
then}
\begin{equation}
\sum_{k=1}^n a_k \sin(kx)>0
\quad{and}
\quad{\sum_{k=0}^n a_k \cos(kx)>0} \quad{(0<x<\pi)}.
\end{equation}

\vspace{0.3cm}
\noindent
In order to prove (1.1) it is enough to consider the special case $a_k=b_k$, where
$$
b_{2k}=b_{2k+1}=\frac{1}{4^k}{2k \choose k}\,\,\,
(k\geq 0)
$$
and to apply summation by parts; see Askey and Steinig \cite{AS}.
In fact, Vietoris proved that
\begin{equation}
S_n(x)>0 \quad\mbox{and}
\quad{T_n(x)>0}
\quad{(n\geq 1; 0<x< \pi)},
\end{equation}
where
$$
S_n(x)=\sum_{k=1}^n b_k \sin(kx)
\quad\mbox{and}
\quad{
T_n(x)=\sum_{k=0}^n b_k\cos(kx)}.
$$
In what follows, we maintain these notations.

In 1974, Askey and Steinig \cite{AS}
offered
a simplified proof of (1.2) and showed that
these inequalities have remarkable applications in the theory of ultraspherical polynomials and that they can be used to find estimates for the location of zeros of trigonometric polynomials.

In the recent past, Vietoris' inequalities received attention from several authors, who offered new conditions on the coefficients $a_k$ such that (1.1) holds;
see Belov \cite{B}, Brown \cite{Br}, 
Brown and Dai \cite{BDW}, Brown and Hewitt \cite{BH}, Brown and Yin \cite{BY}, Koumandos \cite{K}, Mondal and Swaminathan \cite{MS}.
 Interesting historical remarks on these inequalities were given by Askey
\cite{A2}.

Is it possible to replace the lower bound $0$ in (1.2) 
by a positive expression? In 2010, this problem was solved
for the sine polynomial $S_n$. Alzer, Koumandos and Lamprecht
\cite{AKL} proved the following

\vspace{0.3cm}
{\bf{Proposition 2.}}
\emph{The inequalities}
\begin{equation}
S_n(x)>\sum_{k=0}^4 a_k x^k>0
\quad{(a_k\in\mathbf{R}, k=0,...,4)}
\end{equation}
\emph{hold for all $n\geq 1$ and $x\in (0,\pi)$ if and only if}
$$
a_0=0, \quad{a_1=-\pi^2 a_4},
\quad{a_2=3 \pi^2 a_4},
\quad{a_3=-3\pi a_4},
\quad{-1/\pi^3<a_4<0}.
$$
{\em Moreover, in} (1.3) {\em
the biquadratic polynomial cannot be replaced by an algebraic polynomial of degree smaller than $4$.
}

It is natural to ask for a counterpart
of Proposition 2 which holds for the cosine polynomial $T_n$. More precisely, we try to find algebraic polynomials $p$ of smallest degree such that
\begin{equation}
T_n(x)\geq p(x)>0 \quad{(n\geq 1; 0<x<\pi)}.
\end{equation}
It is the aim of this paper to determine
all quadratic polynomials $p$ satisfying
(1.4).

We remark that there is no linear
polynomial $p$  such that (1.4) is valid. Otherwise, 
setting $p(x)=\gamma_0 +\gamma_1 x$
gives 
$$
T_1(x)=1+\cos(x)\geq \gamma_0+\gamma_1 x>0.
$$
We let $x$ tend to $\pi$ and  obtain
$\gamma_0=-\pi \gamma_1$. Thus,
$$
\frac{1+\cos(x)}{\pi-x}\geq -\gamma_1 >0.
$$
But, this contradicts
$$
\lim_{x\to\pi}\frac{1+\cos(x)}{\pi-x}=
0.
$$

In the next section, we collect twelve lemmas. Our main result is presented in Section 3. We conclude the paper with some remarks which
are given in Section 4. Among others,
we provide a new inequality for a sum of
Jacobi polynomials.

The numerical values  in this paper
have been calculated via the computer program MAPLE 13. We point out
that 
in four places we apply the classical Sturm theorem to determine the number of distinct
zeros of an algebraic polynomial in a given interval. Since the Sturm procedure requires lengthy technical computations
we omit the details. However,   those details which we do not include in this paper
are compiled in the supplementary article
\cite{AK}.
Concerning
Sturm's theorem we also refer to
van der Waerden \cite[p. 248 ]{W} and  Kwong \cite{Kw}.

\vspace{0.3cm}
\section{Lemmas}

Here, we collect  lemmas
which play an important role in the proof of
our main result.

\vspace{0.3cm}
{\bf{Lemma 1.}} \emph{We have}
\begin{equation}
\min_{0\leq x<\pi}\frac{T_6(x)}{(x-\pi)^2}
=0.12290... .
\end{equation}

\emph{Proof.}
We define
$$
\eta(x)=10x^6+6x^5-12x^4-\frac{11}{2}x^3+\frac{29}{8}x^2+\frac{11}{8} x+\frac{9}{16}
$$
and
$$
\theta(x)=(\pi-\arccos(x))^2.
$$
Let $c=0.1229$.
First, we show that
\begin{equation}
\eta(x)-c \theta(x)>0 \quad\mbox{for}
\quad{x\in (-1,1]}.
\end{equation}
We distinguish six cases.

Case 1. $-1<x<0$.\\
Let
$$
\omega(x)=\frac{1}{3}x^4+\frac{\pi}{6}x^3+x^2+\pi x +\frac{\pi^2}{4}.
$$
Then we have
\begin{equation}
\omega(x)>0 \quad\mbox{and}
\quad{\omega'(x)>0.}
\end{equation}
Let
$$
\phi(x)=\sqrt{\omega(x)}-\sqrt{\theta(x)}.
$$
Differentiation gives
$$
4(1-x^2)\Bigl(\frac{\omega'(x)}{2\sqrt{\omega(x)}}+\frac{1}{\sqrt{1-x^2}}\Bigr)\omega(x)\phi'(x)=
-\frac{1}{36}x^4(8x+3\pi)(8x^3+3\pi x^2+16x+9\pi)<0,
$$
so that (2.3) leads to
$$
\phi'(x)<0\quad\mbox{and}
\quad{\phi(x)>\phi(0)=0}.
$$
Since
$$
\eta(x)-c\omega(x)>0,
$$
we obtain
$$
\eta(x)-c\theta(x)>c(\omega(x)-\theta(x))>0.
$$

Case 2. $0\leq x\leq 0.3$\\
Using $\eta'(x)>0$ and $\theta'(x)>0$
gives
$$
\eta(x)- c\theta(x)\geq \eta(0)-c\theta(0.3)=0.13... .
$$

Case 3. $0.3\leq x\leq 0.5$.\\
We have $\eta''(x)<0$ and $\theta'(x)>0$. This yields
$$
\eta(x)-c\theta(x)\geq \min (\eta(0.3), \eta(0.5))-c \theta(0.5)=\eta(0.5)-c\theta(0.5)=0.52... .
$$

Case 4. $0.5\leq x\leq 0.65$.\\
Since $\eta'(x)<0$ and $\theta'(x)>0$, we get
$$
\eta(x)-c\theta(x)\geq \eta(0.65)-c\theta(0.65)=0.14... .
$$

Case 5. $0.65\leq x\leq 0.95$.\\
We have $\eta(x)>0$. Let
$$
\lambda(x)=\sqrt{\eta(x)}
-\sqrt{c\theta(x)}
\quad\mbox{and}
\quad{\mu(x)=\eta(x)\eta''(x)-\frac{1}{2}\eta'(x)^2}.
$$
Differentiation leads to
\begin{equation}
2(1-x^2)^3 \Bigl(\frac{\mu(x)}{\eta(x)^{3/2}}+\frac{2\sqrt{c}x}{(1-x^2)^{3/2}}\Bigr)\eta(x)^3\lambda''(x)
=(1-x^2)^{3}\mu(x)^2 -4c x^2 \eta(x)^3=\nu(x), \quad\mbox{say}.
\end{equation}
The function $\mu$ and $\nu$ are polynomials. Applying
Sturm's theorem 
reveals that $\mu$ and
$\nu$ have no zeros on $[0.65,0.95]$.
Since $\mu(3/4)>0$ and $\nu(3/4)>0$, we conclude that both functions are positive
on $[0.65,0.95]$. From (2.4) we find 
 that $\lambda''(x)>0$. Let
$x^*=0.74746$. Then, $\lambda'(x^*)>0$.
This implies
$$
\lambda(x)\geq \lambda(x^*)+(x-x^*)\lambda'(x^*)
\geq \lambda(x^*)+(0.65-x^*)\lambda'(x^*)=0.0000077... .
$$

Case 6. $0.95\leq x\leq 1$.\\
Since $\eta'(x)>0$ and $\theta'(x)>0$, we obtain
$$
\eta(x)-c\theta(x)\geq \eta(0.95)-c\theta(1)=1.43... .
$$
Thus, (2.2) is proved. Let
$$
T^*(x)=\frac{T_6(x)}{(\pi-x)^2}.
$$
We have
$\lim_{x\to\pi} T^*(x)=\infty$ and
$T_6(x)=\eta(\cos(x))$. From (2.2) we obtain
$$
0.1229<T^*(x) \quad{(0\leq x<\pi)}.
$$
Since $T^*(0.725)=0.122907...$, we conclude
that (2.1) is valid.

\vspace{0.3cm}
{\bf{Remark.}}
Numerical computation shows that the minimum value is 
$$
 \alpha  = 0.12290390650...  \label{alp}  
$$ 
attained at the unique point $ x=0. 72656896349... $.

\vspace{0.3cm}
The following lemma is known as 
l'H\^opital's rule for monotonicity. A slightly weaker version can be found in \cite[Proposition 148]{HLP}.

\vspace{0.3cm}
{\bf{Lemma 2.}} \emph{Let $u$ and $v$ be real-valued functions which are continuous
on  $[a,b]$ and differentiable on $(a,b)$.
Furthermore, let $v' \neq 0$ on $(a,b)$. If $u'/v'$ is strictly increasing (resp. decreasing) on $(a,b)$, then 
the functions
$$
x\mapsto \frac{u(x)-u(a)}{v(x)-v(a)}
\quad{and}
\quad{x\mapsto \frac{u(x)-u(b)}{v(x)-v(b)}
}
$$
are strictly increasing (resp. 
decreasing) on $(a,b)$.}

\vspace{0.3cm}
A proof for the next lemma is given by Vietoris in \cite{V}.

\vspace{0.3cm}
{\bf{Lemma 3.}} \emph{Let $n\geq 2$ and $x\in (0,\pi)$.
Then,}
$$
1+\cos(x)-
\frac{1}{4\sin(x/2)}
\Bigl(1+\sin(3x/2)\Bigr)
\leq T_n(x).
$$

\vspace{0.3cm}
{\bf{Lemma 4.}} \emph{If $3\pi/8\leq x<\pi$, then}
\begin{equation}
\frac{123}{1000}
(\pi-x)^2<1+\cos(x)-
\frac{1}{4\sin(x/2)}
\Bigl(1+\sin(3x/2)\Bigr).
\end{equation}

\vspace{0.3cm}
\emph{Proof.}
Let
$$
g(x)=\frac{123}{1000}(\pi-x)^2
\quad\mbox{and}
\quad{h(x)=1+\cos(x)-
\frac{1}{4\sin(x/2)}
\Bigl(1+\sin(3x/2)\Bigr)}.
$$
Then we have 
$$
g'(x)=
\frac{123}{500}(x-\pi)<0
\quad\mbox{and}
\quad{h'(x)=
\frac{\cos(x/2)}{8\sin^2(x/2)}
\Bigl(1-8 \sin^3(x/2)\Bigr)<0.}
$$
Let $3\pi/8\leq r\leq x\leq s\leq 2\pi/3$.
We obtain
$$
h(x)-g(x)\geq h(s)-g(r)=q(r,s),
\quad\mbox{say}.
$$
Since
$$
q(3\pi/8,1.27)=0.0024...,
\quad{q(1.27,1.45)=0.0024...},
\quad{q(1.45,1.7)=0.0008...},
$$
$$
q(1.70,1.95)=0.0072...,
\quad{q(1.95,2\pi/3)=0.0366...},
$$
we conclude that (2.5) holds for $x\in [3\pi/8, 2\pi/3]$. 

Next, we define
$$
f(x)=g(x)-h(x).
$$
Let $y\in (0,\pi/3)$. 
Using
$$
\frac{y}{2}<\frac{\sin(y/2)}{\cos(y/2)}
$$
yields
\begin{eqnarray*}
f'(\pi-y) &=&
\frac{\sin(y/2)}{8 \cos^2(y/2)}
\Bigl( 8 \cos^3(y/2)-1 \Bigr)-
\frac{123 y}{500} \\
  &>&
\frac{\sin(y/2)}{8 \cos^2(y/2)}
\Bigl( 8 \cos^3(y/2)-4\cos(y/2)-1\Bigr) \\[1.2ex] &>& 0.
\end{eqnarray*}
Thus, $y\mapsto f(\pi-y)$ is strictly decreasing on $(0,\pi/3)$, so that
we get
$$
f(\pi-y)< f(\pi)=0.
$$
This proves (2.5) for $x\in (2\pi/3,\pi)$.

\vspace{0.3cm}
{\bf{Lemma 5.}} \emph{For $x\in[0,\pi]$,}
\begin{equation}
x^2<\frac{20000}{99}\Bigl(1-\cos\frac{x}{10}\Bigr).
\end{equation}

\vspace{0.3cm}
\emph{Proof.}
Let
$$
u_0(x)=1-\cos(x/10) \quad\mbox{and}
\quad{v_0(x)=x^2}.
$$
Since
$$
200\,\frac{u_0'(x)}{v_0'(x)}=\frac{\sin(x/10)}{x/10}
$$
is decreasing on $[0,\pi]$, we conclude from Lemma 2 that the function
$$
w(x)=\frac{1-\cos(x/10)}{x^2}
\quad{(0<x\leq \pi)},
\quad{w(0)=\frac{1}{200}}
$$
is also decreasing on $[0,\pi]$. Thus, for $x\in[0,\pi]$,
$$
w(x)\geq w(\pi)=0.004959... >0.00495=\frac{99}{20000}.
$$
This settles (2.6).

\vspace{0.3cm}
{\bf{Lemma 6.}} \emph{Let $a_k$, $\beta_k$ $(k=1,...,n)$, and $\alpha^*$ be real numbers such that}
$$
\sum_{k=1}^j a_k \geq \alpha^*
\quad{for} \quad{j=1,...,n}
\quad{and}
\quad{\beta_1\geq \beta_2 \geq \cdots \geq \beta_n\geq 0}.
$$
\emph{Then,}
$$
\sum_{k=1}^n a_k \beta_k\geq \alpha^* \beta_1.
$$

\vspace{0.3cm}
\emph{Proof.}
Let
$$
A_j=\sum_{k=1}^j a_k \quad\mbox{and}
\quad{\beta_{n+1}=0.}
$$
Summation by parts gives
$$
\sum_{k=1}^n {a_k \beta_k}=\sum_{k=1}^n
A_k(\beta_k -\beta_{k+1})
\geq \sum_{k=1}^n
\alpha^* (\beta_k -\beta_{k+1})=\alpha^* \beta_1.
$$

\vspace{0.3cm}
{\bf{Lemma 7.}} \emph{Let
$$
C_n(x)=\sum_{k=0}^n (-1)^k b_k \cos(kx).
$$
If $2\leq n\leq 21$ $(n\neq 6)$ and $x\in (5\pi/8,\pi)$, then}
\begin{equation}
\frac{820}{33}\Bigl(1-\cos\frac{x}{10}\Bigr)
\leq C_n(x).
\end{equation}

\vspace{0.3cm}
\emph{Proof.}
 We
set $y=x/10$ and
$$
P_n(y)=C_n(10y)-\frac{820}{33}(1-\cos(y)).
$$
Putting $Y=\cos(y)$ reveals that $P_n(y)$ is
 an algebraic polynomial in $Y$.
We denote this polynomial by $P_n^*(Y)$, where
 $Y\in [\cos(\pi/10), \cos(\pi/16)]
=[0.951...,0.980...]$. Applying Sturm's
theorem gives that $P_n^*$ has
no zero on $[0.951, 0.981]$ and satisfies $P_n^*(0.97)>0$. It follows that $P_n$
is positive on $[\pi/16, \pi/10]$. This implies that (2.7) holds.

\vspace{0.3cm}
{\bf{Lemma 8.}} \emph{Let}
\begin{equation}
\Delta(x)=\sum_{k=0}^{21} (-1)^k (b_k -b_{22}) \cos(kx)-\frac{820}{33}\Bigl(1-\cos\frac{x}{10}\Bigr).
\end{equation}
\emph{If \, $5\pi/8\leq x\leq 2.68$, \, then}
\, $\Delta(x)> 0.29$;

\emph{if \, $2.68\leq x\leq 2.83$, \, then}
\, $\Delta(x)> 0.46$;

\emph{if \, $2.83\leq x\leq 2.908$, \, then}
\, $\Delta(x)> 0.64$;

\emph{if \, $2.908\leq x\leq 2.970$, \, then}
\, $\Delta(x)> 0.90$;

\emph{if \, $2.970\leq x\leq 3.021$, \,
then}
\, $\Delta(x)> 1.32$;

\emph{if \, $3.021\leq x\leq 3.051$, \, then}
\, $\Delta(x)> 1.78$.

\vspace{0.3cm}
\emph{Proof.}
Let $5\pi/8\leq x\leq 2.68$. We have
$\cos(\pi/16)=0.980...$ and $\cos(0.268)=0.964...$. 
The function $\Delta -0.29$ is an algebraic
polynomial in $Y=\cos(x/10)$. An application of
Sturm's theorem gives that this function
is positive on $[0.964,0.981]$. This leads
to $\Delta(x)>0.29$ for $x\in [5\pi/8, 2.68]$. Using the same method of proof we obtain that the other estimates for $\Delta(x)$ are also valid.

\vspace{0.3cm}
{\bf{Lemma 9.}} \emph{Let $n\geq 22$,}
\begin{equation}
H_n(x)=\sum_{k=0}^n (-1)^k \cos(kx)
\quad{and}
\quad{D_n(x)=b_{22} H
_{22}(x)+\sum_{k=23}^n (-1)^k b_k\cos(kx).
}
\end{equation}
\emph{If \, $5\pi/8\leq x\leq 2.68$, \, then}
\, $D_n(x)\geq - 0.29$;

\emph{if \, $2.68\leq x\leq 2.83$, \, then}
\, $D_n(x)\geq - 0.46$;

\emph{if \, $2.83\leq x\leq 2.908$, \, then}
\, $D_n(x)\geq - 0.64$;

\emph{if \, $2.908\leq x\leq 2.970$, \, then}
\, $D_n(x)\geq - 0.90$;

\emph{if \, $2.970\leq x\leq 3.021$, \, then}
\, $D_n(x)\geq - 1.32$;

\emph{if \, $3.021\leq x\leq 3.051$, \, then}
\, $D_n(x)\geq - 1.78$.

\vspace{0.3cm}
\emph{Proof.}
We have
$$
H_n(x)
=\frac{1}{2}+(-1)^n\frac{\cos((n+1/2)x)}{2\cos(x/2)}.
$$
Let $5\pi/8\leq x\leq 2.68$. Then we obtain
\begin{equation}
H_n(x)\geq \frac{1}{2}-\frac{1}{2\cos(x/2)}
\geq \frac{1}{2}-\frac{1}{2\cos(1.34)}
=-1.68...>-1.72...=
-\frac{0.29 \cdot 4^{11}}{{22\choose 11}}
=-\frac{0.29}{b_{22}}.
\end{equation}
Using (2.10) and 
$$
b_{22}\geq b_{23} \geq \cdots \geq b_n
$$
we conclude from Lemma 6 that
$D_n(x)\geq -0.29$.
Applying the same method we obtain the other estimates.

\vspace{0.3cm}
As usual, we set
$$
(a)_0=1, \quad{(a)_n=\prod_{k=0}^{n-1}  (a+k)}=\frac{\Gamma(a+n)}{\Gamma(a)}
\quad{(n\geq 1)}.
$$
The following result is due to Koumandos
\cite{K}; see also Koumandos \cite{K1} for background information.

\vspace{0.3cm}
{\bf{Lemma 10.}} \emph{Let  $ 0<\gamma <0. 6915562 $ be given and}
\begin{equation}
d_{2k} = d_{2k+1} = \frac{(\gamma )_k}{k! }
\qquad{(k=0,1,2,...)} .
\end{equation}

\emph{Then, for $n\geq 1$ and $x\in (0,\pi)$,}
$$
\sum_{k=0}^n d_k \cos(kx)>0.
$$

\vspace{0.3cm}
In what follows, we denote by $d_k$ $(k=0,1,2,...)$ the numbers defined in (2.11)
with $\gamma=0.69$.

\vspace{0.3cm}
{\bf{Lemma 11.}} \emph{Let}
\begin{equation}
I(x)= 
\sum_{k=0}^{21}\Bigl(b_k-\frac{b_{22}}{d_{22}} d_k\Bigr)\cos(kx).
\end{equation}
\emph{If $0<x\leq 0.1$, then} $I(x)>1.5$.

\vspace{0.3cm}
\emph{Proof.}
Setting $Y=\cos(x)$ gives that $I-1.5$  is an algebraic polynomial in $Y$. We have $\cos(0.1)=0.995...$. Sturm's theorem reveals that this polynomial is positive
on $[0.995,1]$. It follows that
$I(x)>1.5$ for $x\in[0,0.1]$.

\vspace{0.3cm}
{\bf{Lemma 12.}} \emph{Let $n\geq 22$ and}
\begin{equation}
J_n(x)=\frac{b_{22}}{d_{22}}\sum_{k=0}^{21}
d_k\cos(kx)
+\sum_{k=22}^n
b_k\cos(kx).
\end{equation}
\emph{If $0<x<\pi$, then} $J_n(x)\geq 0$.

\vspace{0.3cm}
\emph{Proof.}
Let $x\in (0,\pi)$.
We set
$$
K_j(x)=\frac{b_{22}}{d_{22}}\sum_{k=0}^j
d_k\cos(kx)\quad{(j=0,1,2,...)}
$$
Then, $K_0(x) \equiv b_{22}/d_{22}$ and from Lemma 10 we obtain
$$
K_j(x)> 0 \quad\mbox{for} \quad{j\geq 1}.
$$
Let
$$
a_k=\frac{b_{22}}{d_{22}} d_k \cos(kx) 
\quad{(k=0,1,...,n)},
\quad{\beta_0 = \cdots =\beta_{21}=1},
\quad{\beta_k=\frac{d_{22} b_k}{b_{22} d_k}} \quad{(k=22,...,n)}.
$$
Since
$$
\beta_0\geq \beta_1 \geq \cdots \geq \beta_n>0,
$$
we conclude from Lemma 6  that
$$
J_n(x)=\sum_{k=0}^n a_k \beta_k  \geq 0.
$$

\vspace{0.3cm}
\section{Main result}

We are now in a position to present
positive lower bounds for the
cosine polynomial $T_n$.

\vspace{0.3cm}
\noindent
{\bf{Theorem.}} \emph{The inequalities
\begin{equation}
T_n(x)\geq c_0 + c_1 x + c_2 x^2 >0
\quad{(c_k\in\mathbf{R}, k=0,1,2)}
\end{equation}
hold for all natural numbers $n$ and real numbers $x\in (0,\pi)$ if and only if}
\begin{equation}
c_0=\pi^2 c_2, 
\quad{c_1=-2\pi c_2,}
\quad{0<c_2\leq \alpha},
\end{equation}
\emph{where}
\begin{equation}
\alpha=\min_{0\leq t<\pi} 
\frac{T_6(t)}{(t-\pi)^2}=0.12290... .
\end{equation}

\vspace{0.3cm}
\emph{Proof.}
We set
$$
Q(x)=c_0+c_1 x+c_2 x^2.
$$
If (3.1) is valid for all $n\geq 1$ and $x\in (0,\pi)$, then we get
$$
T_1(x)=1+\cos(x)\geq Q(x)>0.
$$
We let $x$ tend to $\pi$ and obtain $c_0+c_1\pi +c_2 \pi^2=0$. Thus,
$$
\frac{1+\cos(x)}{x-\pi}\leq\frac{Q(x)}{x-\pi} =c_1+c_2(x+\pi)<0.
$$
Again, we let $x$ tend to $\pi$. This gives
\begin{equation}
c_1=-2\pi c_2 \quad\mbox{and}
\quad{ c_0=-\pi c_1- \pi^2 c_2=\pi^2 c_2.} 
\end{equation}
It follows that
\begin{equation}
Q(x)=c_2(x-\pi)^2 \quad\mbox{with}
\quad{c_2>0}.
\end{equation}
Moreover, from (3.1) (with $n=6$) we obtain
$$
\frac{T_6(x)}{(x-\pi)^2}\geq\frac{Q(x)}{(x-\pi)^2 }=c_2.
$$
Using (3.3) leads to
\begin{equation}
\alpha \geq c_2 .
\end{equation}
From (3.4) - (3.6) we conclude that (3.2) holds.

Next,  we show  that (3.2) and (3.3) lead to
(3.1). If (3.2) and (3.3) are valid, then
$$
0<c_0+c_1 x+c_2 x^2 =c_2(x-\pi)^2
\leq \alpha (x-\pi)^2.
$$
Hence, we have to prove that
\begin{equation} 
\alpha(x-\pi)^2\leq T_n(x)
\quad{(n\geq 1; 0<x<\pi)}.
\end{equation}
Applying Lemma 2 we obtain that the function
$$
F(x)=\frac{T_1(x)}{(x-\pi)^2}=\frac{1+\cos(x)}{(x-\pi)^2}
$$
is strictly increasing on $(0,\pi)$. Thus,
$$
F(x)> F(0)=\frac{2}{\pi^2}=0.202... .
$$
This settles (3.7) for $n=1$. From (3.3) we
conclude that (3.7) is also valid for $n=6$.
In what follows we prove
\begin{equation}
\frac{123}{1000}(\pi-x)^2\leq T_n(x)
\end{equation}
for $n\geq 2$ $(n\neq 6)$ and $x\in (0,\pi)$. With regard to Lemma 3 and Lemma 4
we may assume that $x\in (0,3\pi/8)$. 
We replace in (3.8) $x$ by $\pi-x$. It follows that it is enough to prove
\begin{equation}
\frac{123}{1000}x^2
\leq \sum_{k=0}^n (-1)^k b_k \cos(kx)=C_n(x)
\end{equation}
for $n\geq 2$ $(n\neq 6)$ and $x\in (5\pi/8, \pi)$. 
Using Lemma 5 yields that 
\begin{equation}
\frac{820}{33}\Bigl(1-\cos\frac{x}{10}\Bigr)\leq C_n(x) 
\end{equation}
implies (3.9).
An application of Lemma 7 reveals that (3.10) is valid if $2\leq n\leq 21$ $(n\neq 6)$.

Now, let $n\geq 22$. We have the representation
\begin{equation}
C_n(x)-\frac{820}{33}\Bigl(1-\cos\frac{x}{10}\Bigr)=\Delta(x)+D_n(x),
\end{equation}
with $\Delta$ and $D_n$ as defined in
(2.8) and (2.9), respectively.
Applying Lemma 8 and Lemma 9
reveals that
\begin{equation}
\Delta(x)+D_n(x)\geq 0 \quad\mbox{for}
\quad{x \in [5\pi/8, 3.051].}
\end{equation}
From (3.11) and (3.12) we conclude 
 that
(3.10) is valid for $x\in [5\pi/8,3.051]$.
This implies that (3.8) holds for $x\in [\pi-3.051,3\pi/8]$. Hence, it remains to
prove (3.8) for $x\in (0,\pi-3.051)$.
Since
$$
\frac{123}{1000}(\pi-x)^2<1.22
\quad\mbox{for} \quad{x\in (0,\pi-3.051)}
$$
and $\pi-3.051=0.090...$, it suffices to show that
\begin{equation}
1.22\leq T_n(x) \quad\mbox{for}
\quad{x \in (0, 0.1]}.
\end{equation}
We have
\begin{equation}
T_n(x)=I(x)+J_n(x),
\end{equation}
where $I$ and $J_n$ are defined in (2.12) and (2.13), respectively. Applying Lemma 11 and Lemma 12
we conclude from (3.14) that (3.13) holds. This completes the proof of the Theorem.

\vspace{0.3cm}
\section{Concluding remarks}

(I) If we set $a_0=1$ and 
 $a_k=1/k$ $(k\geq 1)$
in (1.1), then we find 
\begin{equation}
\sum_{k=1}^n \frac{\sin(kx)}{k}>0
\quad\mbox{and}
\quad{1+\sum_{k=1}^n \frac{\cos(kx)}{k}>0}
\quad{(n\geq 1; 0<x<\pi)}.
\end{equation}
The first inequality  is the famous
Fej\'er-Jackson inequality, which was conjectured by Fej\'er in 1910 and proved
one year later by Jackson \cite{J}. Its analogue for the cosine sum was published by Young \cite{Y} in 1913. Both inequalities  motivated the research of many
authors, who  presented numerous refinements, extensions, and variants
of (4.1). We refer to Askey \cite{A1}, Askey and Gasper \cite{AG},
Milovanovi\'c, Mitrinovi\'c, and Rassias \cite[chapter 4]{MMR} and the references cited therein.

Applying our Theorem
 we obtain an improvement of Young's inequality:
$$
1+\sum_{k=1}^n \frac{\cos(kx)}{k}
\geq \alpha (\pi-x)^2 \quad{(n\geq 1; 0<x<\pi)},
$$
{where $\alpha$ is given in} (3.3)

\vspace{0.2cm}
(II)  Askey and Steinig \cite{AS}
used Proposition 1 to prove the following
interesting result.

\vspace{0.3cm}
{\bf{Proposition 3.}} \emph{Let $\gamma_k$
$(k=0,1,...,n)$ be positive real numbers such that}
\begin{equation}
2k \gamma_k \leq (2k-1) \gamma_{k-1}
\quad{(k\geq 1)}.
\end{equation}
\emph{Then, for $n\geq 0$ and $t\in (0,2\pi)$,}
$$
\sum_{k=0}^n \gamma_k \sin \bigl( (k+1/4)t\bigr)>0
\quad{and} 
\quad{\sum_{k=0}^n \gamma_k \cos\bigl( (k+1/4)t\bigr)>0.}
$$

An application of the Theorem leads to a refinement of the second inequality:
\begin{equation}
\sum_{k=0}^n \gamma_k \cos \bigl( (k+1/4)t\bigr)\geq  \frac{\alpha \gamma_0(2\pi-t)^2}{8 \cos(t/4)}
\quad{(n\geq 0; 0<t<2\pi)}.
\end{equation}

In order to prove (4.3) we set
$$
\gamma_k^*=\frac{1}{4^k}{2k\choose k}
=\frac{1}{k! }\Bigl(\frac{1}{2}\Bigr)_k
\quad{(k\geq 0)}.
$$
Then,
\begin{equation}
2\cos(t/4)\sum_{k=0}^j \gamma_k^* \cos \bigl( (k+1/4)t\bigr)=T_{2j+1}(t/2)
\geq \alpha (\pi-t/2)^2
\quad\mbox{for}
\quad{j=0,1,...,n}.
\end{equation}
We define
$$
\beta_k=\frac{\gamma_k}{\gamma_k^*}
\quad{(k\geq 0)}
$$
and apply (4.2). This yields
$$
\beta_0\geq \beta_1 \geq \cdots \geq \beta_n>0.
$$
Using (4.4) and Lemma 6 leads to
$$
\sum_{k=0}^n \gamma_k \cos \bigl( (k+1/4)t\bigr)=
\sum_{k=0}^n \beta_k\gamma_k^* \cos\bigl( (k+1/4)t\bigr)
\geq \beta_0 \frac{\alpha(\pi-t/2)^2}{2\cos(t/4)}.
$$
Since $\beta_0=\gamma_0$, we get (4.3).

\vspace{0.2cm}
(III) The classical Jacobi polynomials
$P_m^{(a,b)}(z)$ are given by
$$
P_m^{(a,b)}(z)=
\frac{(a+1)_m}{m! }\sum_{k=0}^m
\frac{(-m)_k (m+a+b+1)_k}{k! (a+1)_k}\Bigl(\frac{1-z}{2}\Bigr)^k.
$$
A collection of the main properties of these functions
can be found, for instance, in 
 \cite[chapter 1.2.7]{MMR}.
An application of (4.4) (with $t=4x, j=n)$ and the
identity
$$
\frac{P_m^{(-1/2,-1/2)}(\cos(x))}{P_m^{(-1/2,-1/2)}(1)}=\cos(m x)
$$
(with $m=4k+1)$ yields
$$
\sum_{k=0}^n \frac{1}{k! } \Bigl(\frac{1}{2}\Bigr)_k
\frac{P_{4k+1}^{(-1/2,-1/2)}(\cos(x))}{P_{4k+1}^{(-1/2,-1/2)}(1)}\geq 
\frac{\alpha (\pi-2 x)^2}{2\cos(x)}
\quad{(n\geq 0; 0< x<\pi/2)}.
$$
For related inequalities we refer to
Askey \cite{As},  Askey and Gasper
\cite{AG1}, \cite{AG} and the references therein.

\vspace{0.7cm}
{\bf{Acknowledgement.}} We thank the referees for helpful comments.

\end{document}